\documentclass[runningheads]{llncs}
\usepackage[T1]{fontenc}
\usepackage{graphicx}
\usepackage{amssymb, latexsym, amsmath, amsfonts, epsfig}
\newcommand{\bfs}[1]{{\boldsymbol #1}}

\begin{document}
\title{Isogeometric Analysis of Bound States of a Quantum Three-Body Problem in 1D
}
\titlerunning{IGA of Three-Body Problems}
%
\author{Quanling Deng\inst{1} }
%
\authorrunning{Q. Deng}

%
\institute{School of Computing, Australian National University, Canberra, ACT 2601, Australia. \email{Quanling.Deng@anu.edu.au} }
%
\maketitle              
\begin{abstract}
In this paper, we initiate the study of isogeometric analysis (IGA) of a quantum three-body problem that has been well-known to be difficult to solve. 
In the IGA setting, we represent the wavefunctions by linear combinations of B-spline basis functions and solve the problem as a matrix eigenvalue problem. 
The eigenvalue gives the eigenstate energy while the eigenvector gives the coefficients of the B-splines that lead to the eigenstate. 
The major difficulty of isogeometric or other finite-element-method-based analyses lies in the lack of boundary conditions and a large number of degrees of freedom required for accuracy.  
For a typical many-body problem with attractive interaction, there are bound and scattering states where bound states have negative eigenvalues.
We focus on bound states and start with the analysis for a two-body problem. 
We demonstrate through various numerical experiments that IGA provides a promising technique to solve the three-body problem. 

\end{abstract}
\section{Introduction} \label{sec:intr}

While there are still unsolved questions in the classical three-body problem \cite{breen2020newton,letellier2019chaos}, 
the quantum mechanical three-body problem also has unanswered questions \cite{drut2018quantum,sukhareva2018validity,guo2018numerical}. 
Accurate and rigorous solutions are highly desirable both for answering these open questions as well as for studies of three-body correlations within various many-body systems.
The two-body problem is generally considered as ``solved" due to the momentum conservation that leads to a second-order ordinary differential equation which can be solved effectively. 
For three-body problem, the space is six-dimensional in the center of mass system.
The total angular momentum conservation leads to three coupled second-order nonlinear differential equations 
in classical mechanics \cite{nielsen2001three}.  
In quantum mechanics, we have a more complicated system that admits no analytic solution in general 
and existing numerical methods are not  satisfactory in the sense of robustness, efficiency, and reliability.

In quantum chemistry and molecular physics, the Born--Oppenheimer (BO) approximation has been the most well-known and widely-used mathematical approximation since the early days of quantum mechanics \cite{baer2006beyond,panati2007time,cederbaum2008born,scherrer2017mass}. 
The method is based on the assumption that the nuclei are much heavier than the electrons which consequently leads to that the wave functions of atomic nuclei and electrons in a molecule can be treated separately.
For instance, BO was used in \cite{fonseca1979efimov,efremov2009efimov} to study the Efimov effect in few-body systems.
BO was adopted in \cite{happ2019universality,happ2022universality} recently to  establish the universality in a one-dimensional three-body system.
A more efficient numerical method based on BO was developed using the tensor-product structure \cite{thies2021exploiting}. 
BO has been the standard method to describe the interaction between electrons and nuclei but it can fail whenever the assumption fails, for example, in graphene \cite{pisana2007breakdown}. 
When the mass ratios of the interacting bodies are close to one, the assumption fails and BO is generally inaccurate. 
Other methods such as the pseudospectral method based on Fourier analysis \cite{boyd2001chebyshev}  
and Skorniakov and Ter-Martirosian (STM) method  based on exact integral equations \cite{skorniakov1957three} have been developed to obtain the three-body
bound states for arbitrary mass ratios.

In this paper, we develop a general numerical method to solve one-dimensional quantum two- and three-body problems with arbitrary mass ratios 
and any interaction potentials that lead to bound states.
With this goal in mind, we initiate the study of finite element analysis (FEA) based methods to find the bound states of three-body systems. 
In particular, we adopt the more advanced method isogeometric analysis (IGA) for this purpose. 
IGA, first developed in \cite{hughes2005isogeometric,cottrell2009isogeometric}, has been widely-used as a numerical analysis tool for various simulations that are governed by partial differential equations (PDEs).
IGA adopts the framework of classic Galerkin FEA and uses B-splines or non-uniform rational basis splines (NURBS) instead of the Lagrange polynomials as its basis functions. 
These basis functions have higher-order continuity (smoothness) which consequently improves the accuracy of the FEA numerical approximations.
The work  \cite{cottrell2006isogeometric} applied IGA to study a structural vibration problem that is modeled as a Laplacian eigenvalue problem. 
It has been shown that IGA improved the accuracy of the spectral approximation significantly compared with FEA \cite{hughes2008duality}. 
Further advantages of IGA over FEA on spectral accuracy have been studied in \cite{hughes2014finite,puzyrev2017dispersion}.
With the advantages in mind, we adopt IGA to solve the quantum three-body problem as a second-order differential eigenvalue problem.

The rest of this paper is organized as follows.  
Section~\ref{sec:ps} presents the two- and three-body problems under consideration.
We then unify these two problems as a single differential eigenvalue problem in one or two dimensions where 1D refers to the two-body problem and 2D refers to the three-body problem. 
We show an example of solutions to a two-body problem for both bound and scattering states, which serves as a motivation of the proposed method that solves only the bound states over an approximate finite domain. 
We then present the IGA discretization method in Section \ref{sec:iga} to solve the unified problem for the bound states. 
Section \ref{sec:num} collects and discusses various numerical tests to demonstrate the performance of the proposed method.
We also perform the numerical study of the impact of domain size on the approximation accuracy of the bound states.
Concluding remarks are presented in Section~\ref{sec:conclusion}.

\section{The two- and three-body problems} \label{sec:ps}
In this section, we first present the heavy-light two-body and heavy-heavy-light three-body problems 
that are modeled as the dimensionless stationary Schr\"{o}dinger equations recently studied in \cite{happ2019universality,happ2022universality}. 
We then generalize the problems for any mass ratios and unify them as a single differential eigenvalue problem. 
A numerical example is followed to show the bound and scattering states of a two-body problem. 
The shape of bound states gives a motivation to pose the differential eigenvalue problem on a finite domain with a size to be specified depending on the differential operator and accuracy tolerance. 

The heavy-light quantum two-body system with an attractive interaction via a potential of finite range,
after eliminating the center-of-mass motion, is modeled as a dimensionless
stationary Schr\"{o}dinger equation
\begin{equation} \label{eq:2bp}
\Big[ - \frac{1}{2} \frac{\partial^2}{\partial x^2} - v(x) \Big] \psi^{(2)}  = E^{(2)} \psi^{(2)},
\end{equation}
where $E^{(2)}$ is the binding energy and $\psi^{(2)}$ is the two-body wave function. 
The corresponding three-body system is modeled as
\begin{equation} \label{eq:3bp}
\Big[ - \frac{\alpha_x}{2} \frac{\partial^2}{\partial x^2} -  \frac{\alpha_y}{2} \frac{\partial^2}{\partial y^2} - v(x + y/2) - v(x-y/2) \Big] \psi  = E \psi,
\end{equation}
where $E$ is the eigenenergy and $\psi = \psi(x,y)$ is the three-body wave function describing the relative motions. 
The coefficients
\begin{equation}
\alpha_x = \frac{1/2 + m_h/m_l}{1 + m_h/m_l}, \qquad  \alpha_y = \frac{2}{1 + m_h/m_l},
\end{equation}
where $m_h$ denotes the mass of two heavy particles and $m_l$ denotes the mass of the light particle. 
The potential 
\begin{equation}
v(\xi) = \beta f(\xi) 
\end{equation}
with $\beta>0$ denoting a magnitude and $f$ denoting the shape of the interaction potential. 
We assume that the $f$ is symmetric and describes a short-range interaction, that is, $|\xi |^2 f( |\xi | ) \to 0$ as $|\xi | \to \infty$.

\subsection{The unified problem} \label{sec:up}
The two-body problem \eqref{eq:2bp} is posed on an infinite domain $\Omega = \mathbb{R}$ while the three-body problem \eqref{eq:3bp} is posed on $\Omega = \mathbb{R}^2$.
Mathematically, problems \eqref{eq:2bp} and \eqref{eq:3bp} are differential eigenvalue problems where the differential operator is a Hamiltonian. 
Moreover, we observe that \eqref{eq:2bp} is of one variable while \eqref{eq:3bp}  is of two variables which can be regarded as
1D and 2D spatial variables, respectively.
With this in mind, we unify problems \eqref{eq:2bp} and \eqref{eq:3bp} to obtain a differential eigenvalue problem
\begin{equation} \label{eq:pde0}
\begin{aligned}
- \nabla \cdot (\kappa \nabla u) - \gamma u& =  \lambda u \quad &&  \forall \ \bfs{x} \in \Omega, \\
\end{aligned}
\end{equation}
where $\nabla$ is the gradient operator, $\nabla \cdot$ is the divergence operator. $\lambda = E^{(2)}, \gamma = v(x), \kappa = \frac12$ in 1D while $\lambda = E, \gamma = v(x + y/2) + v(x-y/2)$ and $\kappa = (\frac{\alpha_x}{2}, 0; 0, \frac{\alpha_y}{2})$ being a diagonal matrix in 2D. Herein, $u$ denotes an eigenstate.

From now on, we focus on the unified problem \eqref{eq:pde0}. 
There are three major difficulties in solving this problem using a Galerkin FEA-based discretization method.
\begin{itemize}
\item (a) The attractive interaction may lead to negative eigenvalues. Consequently, the discretization of the differential operator $\mathcal{L} = - \nabla \cdot (\kappa \nabla) - \gamma$ leads to a stiffness matrix that is not necessarily positive-definite. This in return brings a potential issue when solving the resulting linear algebra problem. 
\item (b) The domain $\Omega$ is infinite. This makes it impossible to discretize the domain with a finite number of elements with each element being of finite size. 
\item (c) There are no boundary conditions provided. A Galerkin FEA-based discretization method requires setting appropriate boundary conditions for the resulting linear algebra system to be non-singular. 
\end{itemize}

For (a), an eigenvalue shift will resolve the issue. That is, we rewrite $- \nabla \cdot (\kappa \nabla u) - \gamma u =  \lambda u$ by adding a positive scale to obtain $- \nabla \cdot (\kappa \nabla u) - (\gamma -\gamma_0) u =  (\lambda + \gamma_0) u$ where $\gamma_0 > 0$ is a constant such that $\gamma -\gamma_0<0$ for all $\boldsymbol{x} \in \Omega$. 
With a slight abuse of notation, the problem \eqref{eq:pde0} can be rewritten as 
\begin{equation} \label{eq:pde}
\begin{aligned}
- \nabla \cdot (\kappa \nabla u) + \gamma u& =  \lambda u \quad &&  \forall \ \bfs{x} \in \Omega.
\end{aligned}
\end{equation}
To overcome the difficulties (b) and (c), we first present an example of a solution to the two-body problem in the next subsection.

\subsection{A solution example of the two-body problem} \label{sec:ex}
For attractive interaction, the eigenenergies for certain eigenstates can be negative. 
For a potential vanishing at $\pm \infty$, a negative eigenvalue implies a bound state while a positive eigenvalue implies a scattering state \cite{griffiths2018introduction}. 

\begin{figure}[h!]
\hspace{-0.5cm}
\includegraphics[height=5.5cm]{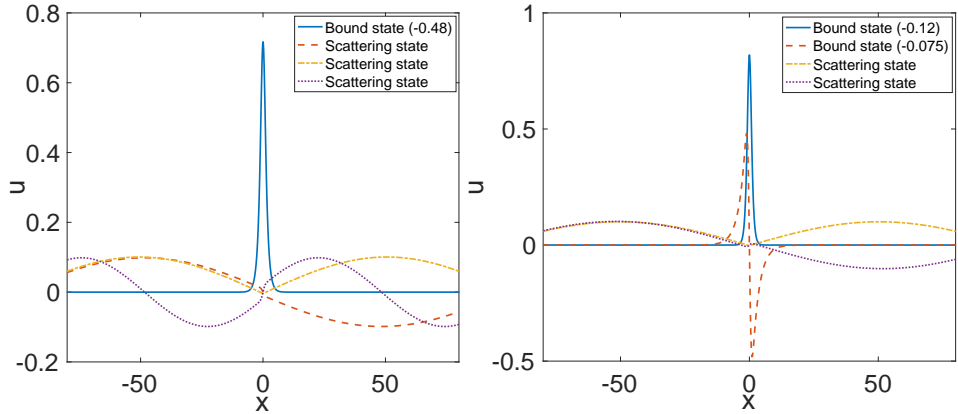}
\vspace{-0.5cm}
\caption{Bound and scattering states of two-body problems.}
\label{fig:2bp}
\end{figure}

Figure \ref{fig:2bp} shows an example of state solutions to the two-body problem \eqref{eq:2bp} 
with $\kappa=1/2$, potential $f(\xi) = e^{-\xi^2}$ and $\beta=1$ for the left plot while $\beta=2$ for the right plot. 
The bound states eigenenergies are marked in the figure. 
Herein, we apply IGA with 5000 elements and $C^6$ septic B-spline basis functions. 
We present the details of the IGA method in the next section.
We observe that there is one bound state for $\beta=1$ and two bound states for $\beta=2$. 
All other states are scattering states. 
When $x\to \pm \infty$, the wavefunctions go to zeros exponentially fast for bound states while they do not go to zeros for scattering states. A theoretical explanation can be found in \cite{agmon2014lectures}.

In this paper, our goal is to find the eigenenergies and eigenstates for bound states. 
In the case of $\beta=1$, the bound state solution decays to zero approximately at $x=\pm 10$ with an error of $5.9 \times 10^{-5}$. 
We observe similar behaviour for the two bound states of the case $\beta=2$.
This decaying behaviour provides an insight to overcome the difficulties (b) and (c) listed in section \ref{sec:up}.
The idea is that for a given tolerance $\epsilon >0$, we propose to solve the problem \eqref{eq:pde} on a finite domain
$\Omega_\epsilon = [-\boldsymbol{x}_\epsilon, \boldsymbol{x}_\epsilon]^d, d=1,2$ with homogeneous boundary condition 
\begin{equation} \label{eq:bc}
u = 0, \qquad \forall \ \bfs{x} \in \partial \Omega_\epsilon.
\end{equation}

\begin{remark}
For smaller error tolerance, one expects to apply a larger finite domain. 
A detailed study is presented in section \ref{sec:ds}.
We also point out that a simple transformation such as $\hat x = \tanh(x)$ that transfers the infinite domain $\Omega$ to a finite domain $\hat\Omega = (-1, 1)^d, d=1,2$ can not avoid the difficulties of (b) and (c) listed in section \ref{sec:up}.
\end{remark}

\section{Isogeometric analysis} \label{sec:iga}

In this section, we present the IGA method for the unified problem \eqref{eq:pde} on $\Omega_\epsilon$ supplied with the boundary condition \eqref{eq:bc}. 
We also give an a priori error estimate for the bound states and their eigenenergies. 

\subsection{Continuous level}
Let $\Omega_\epsilon = [-\boldsymbol{x}_\epsilon, \boldsymbol{x}_\epsilon]^d \subset  \mathbb{R}^d, d=1,2$ be a bounded domain 
with Lipschitz boundary $\partial \Omega_\epsilon$. 
We adopt the standard notation for the Hilbert and Sobolev spaces. 
For a measurable subset $S\subseteq \Omega_\epsilon$, 
we denote by $(\cdot,\cdot)_S$ and $\| \cdot \|_S$ the $L^2$-inner product and its norm, respectively. 
We omit the subscripts when clear.
For an integer $m\ge1$,
we denote the $H^m$-norm and $H^m$-seminorm as $\| \cdot \|_{H^m(S)}$ and 
$| \cdot |_{H^m(S)}$, respectively. In particular, we denote by $H^1_0(\Omega_\epsilon)$ the Sobolev space 
with functions in $H^1(\Omega_\epsilon)$ that are vanishing on the boundaries.

The variational formulation of \eqref{eq:pde} is to find eigenvalue $\lambda \in \mathbb{R}^{+}$ and 
eigenfuction $u \in H^1_0(\Omega_\epsilon)$ with $\|u\|_{\Omega_\epsilon}=1$ such that 
\begin{equation} \label{eq:vf}
a(w, u) =  \lambda b(w, u), \quad \forall \ w \in H^1_0(\Omega_\epsilon), 
\end{equation}
where the bilinear forms are defined as for $v,w \in H^1_0(\Omega_\epsilon)$
\begin{equation}
a(v,w) := (\kappa \nabla v, \nabla w)_{\Omega_\epsilon} + (\gamma v, w)_{\Omega_\epsilon},
\qquad
b(v,w) := (v,w)_{\Omega_\epsilon}.
\end{equation}

The eigenvalue problem \eqref{eq:vf} with $\gamma = \gamma(\bfs{x}) >0, \forall \bfs{x} \in \Omega_\epsilon,$ has a countable set of positive eigenvalues 
(see, for example, \cite[Sec. 9.8]{Brezis:11}) 
\begin{equation*}
0 < \lambda_1 < \lambda_2 \leq \lambda_3 \leq \cdots
\end{equation*}
with an associated set of orthonormal eigenfunctions $\{ u_j\}_{j=1}^\infty$. Thus, there holds 
$
(u_j, u_k) = \delta_{jk}, 
$
where $\delta_{jk} =1$ is the Kronecker delta. 
As a consequence, the eigenfunctions are also orthogonal in the energy inner-product as there holds 
$
a(u_j, u_k) =  \lambda_j b(u_j, u_k) = \lambda_j \delta_{jk}.
$

\subsection{IGA discretized level}
At the discretized level, we first discretize the domain $\Omega_\epsilon$ with a uniform tensor-product mesh. 
We denote a general element as $\tau$ and its collection as $\mathcal{T}_h$ such that 
$\overline \Omega_\epsilon = \cup_{\tau\in \mathcal{T}_h}  \tau$. 
Let $h = \max_{\tau \in \mathcal{T}_h} \text{diameter}(\tau)$. 
In the IGA setting, for simplicity, we use the B-splines. 
The B-spline basis functions in 1D are given as the Cox-de Boor recursion formula; 
we refer to \cite{de1978practical,piegl2012nurbs} for details.
Let $X = \{x_0, x_1, \cdots, x_m \}$ be a knot vector with a nondecreasing sequential knots $x_j$.  
The $j$-th B-spline basis function of degree $p$, denoted as $\phi^j_p(x)$, is defined recursively as 
\begin{equation} \label{eq:Bspline}
\begin{aligned}
\phi^j_0(x) & = 
\begin{cases}
1, \quad \text{if} \ x_j \le x < x_{j+1}, \\
0, \quad \text{otherwise}, \\
\end{cases} \\ 
\phi^j_p(x) & = \frac{x - x_j}{x_{j+p} - x_j} \phi^j_{p-1}(x) + \frac{x_{j+p+1} - x}{x_{j+p+1} - x_{j+1}} \phi^{j+1}_{p-1}(x).
\end{aligned}
\end{equation}
A tensor-product of these 1D B-splines produces the B-spline basis functions in multiple dimensions. 
We define the multi-dimensional approximation space as $V^h_p \subset H^1_0(\Omega_\epsilon)$ with:
\begin{equation*}
V^h_p = \text{span} \{ \phi^j_p \}_{j=1}^{N_h} = 
\begin{cases}
 \text{span} \{ \phi^{j_x}_{p_x}(x) \}_{j_x=1}^{N_x}, & \text{in 1D}, \\
 \text{span} \{ \phi^{j_x}_{p_x}(x) \phi^{j_y}_{p_y}(y) \}_{j_x, j_y=1}^{N_x, N_y}, & \text{in 2D}, \\
\end{cases}
\end{equation*}
where $p_x, p_y$ specify the approximation order in each dimension. 
$N_x, N_y$ is the total number of basis functions in each dimension 
and $N_h$ is the total number of degrees of freedom.
The isogeometric analysis of \eqref{eq:pde} in variational formulation seeks $\lambda^h \in \mathbb{R}$ and $u^h \in V^h_p$ with $\| u^h \|_{\Omega_\epsilon} = 1$ such that 
\begin{equation} \label{eq:vfh}
a(w^h, u^h) =  \lambda^h b(w^h, u^h), \quad \forall \ w^h \in V^h_p.
\end{equation}

\subsection{Algebraic level}
At the algebraic level, we approximate the eigenfunctions as a linear combination of the B-spline basis functions, i.e.,
$$
u^h = \sum_{j=1}^{N_h} \nu_j \phi^j_p,
$$
where $\nu_j,j=1,\cdots, N_h$ are the coefficients.  
We then substitute all the B-spline basis functions for $w^h$ in \eqref{eq:vfh}. 
This leads to the generalized matrix eigenvalue problem
\begin{equation} \label{eq:mevp}
\mathbf{K} \mathbf{U} = \lambda^h \mathbf{M} \mathbf{U},
\end{equation}
where $\mathbf{K}_{kl} =  a(\phi_p^k, \phi_p^l), \mathbf{M}_{kl} = b(\phi_p^k, \phi_p^l),$ 
and $\mathbf{U}$ is the corresponding representation of the eigenvector as the coefficients of the B-spline basis functions. 
The homogeneous Dirichlet boundary condition \eqref{eq:bc} can be set by removing the rows and columns corresponding to the degrees of freedom associated with the nodes at the boundary. 
This matrix eigenvalue problem is to be solved in a computing program.

\subsection{A priori error estimates} \label{sec:ee}
IGA is a Galerkin finite element discretization method. 
On a rectangular domain with tensor-product grids, 
the only difference of IGA from the classical FEA is the basis functions. 
FEA adopts $C^0$ polynomials as basis function while IGA adopts $C^k, k\ge 1$ polynomials. 
We observe that $C^0$ is a larger space, i.e., $C^k(\Omega_\epsilon) \subset C^0(\Omega_\epsilon)$.
In general, for an a priori error estimate that is established in the Galerkin FEA framework, 
the estimate also holds for IGA. 
Thus, we expect optimal convergence rates for the eigenvalues and eigenfunctions as in FEA \cite{ciarlet1978finite,babuvska1991eigenvalue,Ern_Guermond_FEs_II_2020}. 
We present the following estimate without a theoretical proof. 
Instead, we show numerical validation in Section \ref{sec:num}.

Given the mesh configuration and IGA setting described above, let $(\lambda_j, u_j) \in \mathbb{R}^+\times H^1_0(\Omega)$ solve \eqref{eq:vf} for  bound states
and let $(\lambda_j^h,  u_j^h) \in \mathbb{R}^+\times V^h_p$ solve \eqref{eq:vfh} for  bound states with the normalizations $\| u_j \|_{\Omega_\epsilon} = 1$ and $ \| u_j^h \|_{\Omega_\epsilon} = 1$. 
Assuming elliptic regularity on the operator $\mathcal{L} = - \nabla \cdot (\kappa \nabla) + \gamma$ 
and high-order smoothness of the eigenfunctions $u_j$ on $\Omega_\epsilon$, there holds:
\begin{equation} \label{eq:ee}
\big| \lambda_j^h - \lambda_j \big|  \le Ch^{2p}, \qquad | u_j - u_j^h |_{H^1(\Omega)} \le Ch^p, 
\end{equation}
where $C$ is a positive constant independent of the mesh-size $h$. 
We remark that these estimates only hold for bound states and do not necessarily hold for scattering states.

\section{Numerical experiments} \label{sec:num}

In this section, we present various numerical examples to demonstrate the performance of IGA. 
We first show the IGA approximation optimal convergence accuracy with a domain $\Omega_\epsilon$ of large size.  
Then we study the impact of the domain size on accuracy and 
give an approximate formula that determines the size of the domain given a certain accuracy tolerance.

\subsection{IGA discretization accuracy}

We focus on the two- and three-body problems with a potential with polynomial decay 
\begin{equation} \label{eq:pd}
f(\xi) = \frac{1}{(1 + \xi^2 )^3}
\end{equation} 
of the cube of a Lorentzian and one with exponential decay 
\begin{equation} \label{eq:ed}
f(\xi) = e^{-\xi^2}
\end{equation}
of a Gaussian. 
For these potentials, finding the exact analytical solutions is impossible. 
For the purpose of characterizing the errors, we use, as a reference solution to the exact one,
 the solution of IGA with a septic $C^6$ B-spline basis functions and a fine mesh. 
We focus on the eigenvalue error that is defined as
\begin{equation}
e_j = |  \lambda^h_j - \hat \lambda_j |,
\end{equation}
where $\lambda^h_j$ is an IGA eigenvalue and $\hat \lambda_j$ is a reference eigenvalue that is of high accuracy approximating the exact one $\lambda_j$.

\begin{figure}[h!]
\hspace{-0.5cm}
\includegraphics[height=6cm]{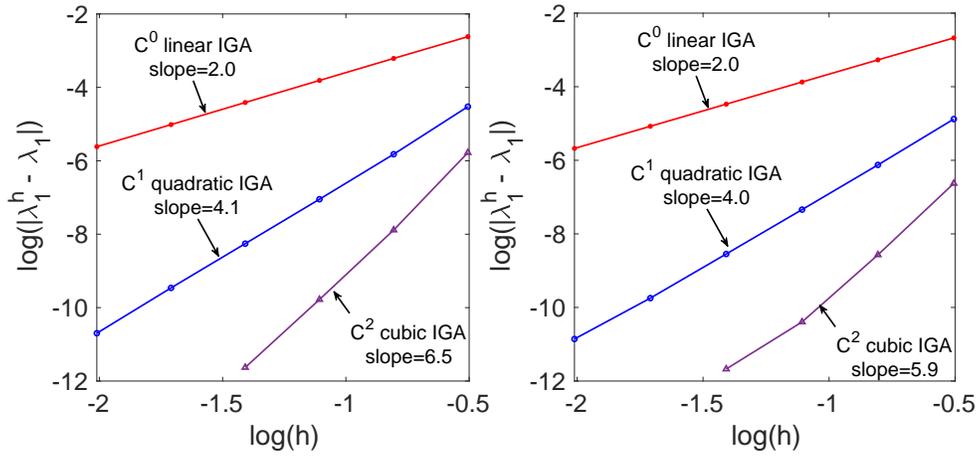}
\caption{Eigenvalue error convergence rates of IGA for the two-body problem with $\beta=1$ on a domain $\Omega_\epsilon = [-20, 20].$ The potential has polynomial decay  \eqref{eq:pd} for the left plot while exponential decay \eqref{eq:ed} for the right plot.  }
\label{fig:beta1err}
\end{figure}

Figure \ref{fig:beta1err} shows the eigenvalue error convergence rates for IGA of the two-body problem with $\kappa=\frac12$ in 1D. 
We consider $C^0$ linear, $C^1$ quadratic, and $C^2$ cubic IGA elements. 
We study the problem with both potentials  \eqref{eq:pd}  and \eqref{eq:ed} and a fixed magnitude $\beta=1$.
For both potentials \eqref{eq:pd}  and \eqref{eq:ed}, there is one bound state. 
The state reference eigenvalue is $\hat \lambda_1 = -0.31658012845$ for \eqref{eq:pd}  and $\hat \lambda_1 = -0.47738997738$ for \eqref{eq:ed}, respectively.
We solve the problem for the bound state using $C^6$ septic IGA with 5000 uniform elements over the domain $\Omega_\epsilon = [-20,20]$.
We observe optimal error convergence rates in all the scenarios. 
This confirms the theoretical prediction \eqref{eq:ee} in section \ref{sec:ee}. 

\begin{figure}[h!]
\hspace{-0.7cm}
\includegraphics[height=5.5cm]{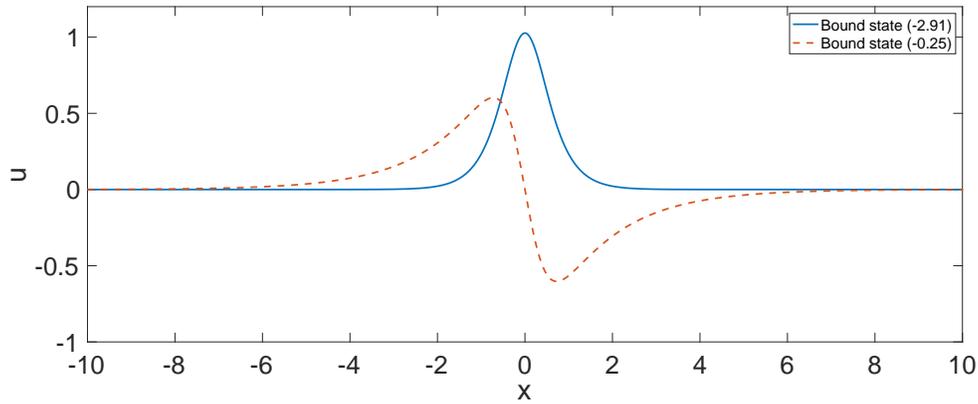}
\caption{The bound states of the two-body problem with $\beta=5$ and a potential of polynomial decay \eqref{eq:pd}.  }
\label{fig:2bpb5}
\end{figure}

Now we consider a case where there are two bound states in the two-body problem. Let $\beta=5$ and we apply a potential of polynomial decay \eqref{eq:pd}. 
Figure \ref{fig:2bpb5} shows the two bound states solutions and their eigenenergies are $-2.9149185630$ and $-0.25417134380$. 
Herein, the numerical eigenstates are computed using $C^6$ septic IGA with 5000 uniform elements over the domain $\Omega_\epsilon = [-20,20]$. The plot shows the eigenstate over $[-10,10]$ for better focus while the problem is solved over the larger domain $\Omega_\epsilon = [-20,20]$ for high accuracy.

\begin{figure}[h!]
\hspace{-0.5cm}
\includegraphics[height=6cm]{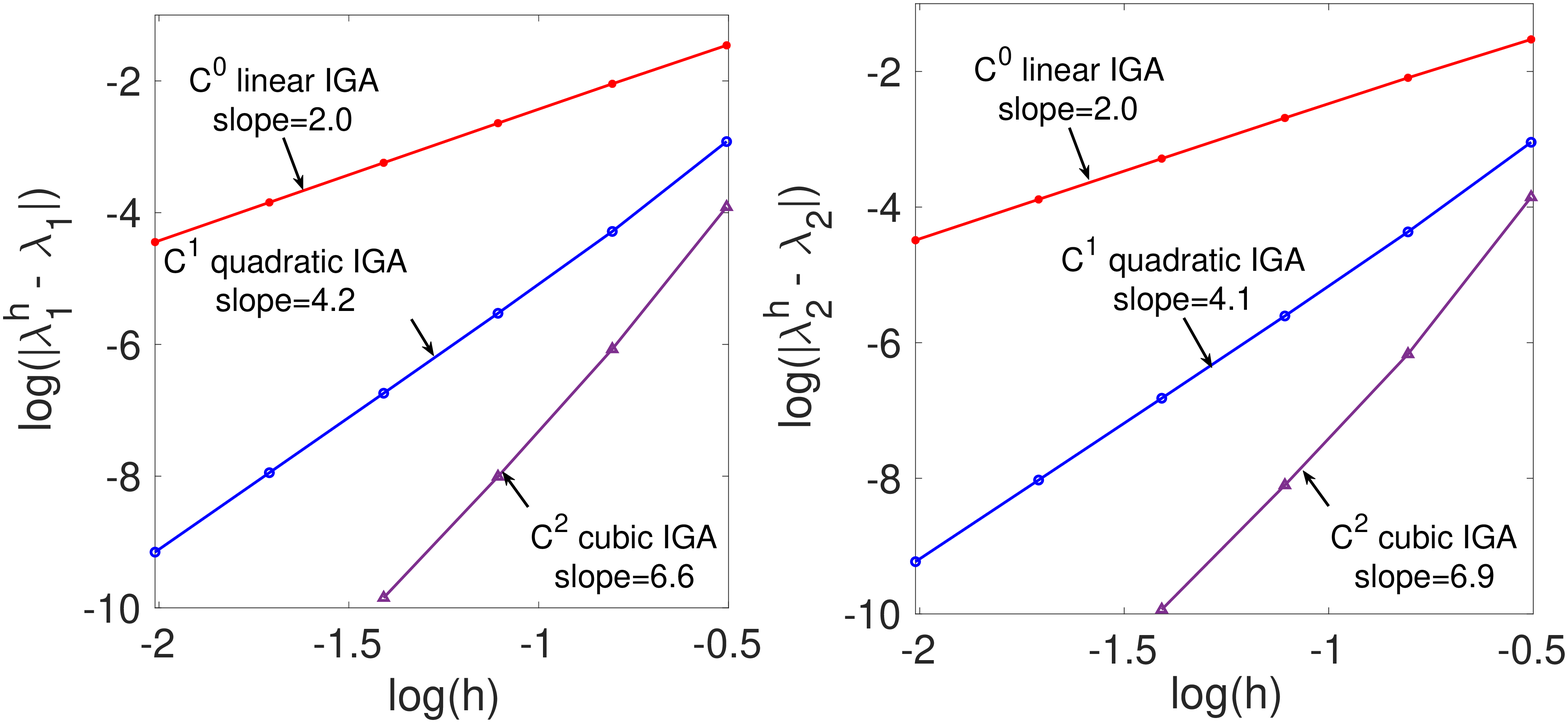}
\vspace{-0.5cm}
\caption{
Eigenvalue error convergence rates of IGA for the two-body problem with $\beta=5$ and a  polynomially decaying potential \eqref{eq:pd} on a domain $\Omega_\epsilon = [-20, 20].$  The left plot shows the eigenenergy errors of first bound state while the right plot shows these of the second bound state.}
\label{fig:beta5err}
\end{figure}

Figure \ref{fig:beta5err} shows the eigenenergy error convergence rates of the problem described above (also shown in Figure \ref{fig:2bpb5}). Again, we observe optimal error convergence rates that verify the theoretical prediction. 
Moreover, for IGA with higher-order elements, the eigenvalue errors reach small errors faster with coarser meshes. 
This validates that the reference solutions obtained by using $C^6$ septic IGA with 5000 elements are of high accuracy and can be used as highly accurate approximations to the exact solutions.

\subsection{A study on domain size} \label{sec:ds}

For bound states of the two-body problems as discussed in section \ref{sec:ex}, the state values approach zero exponentially fast. 
The IGA discretization requires a finite domain with homogeneous boundary condition \eqref{eq:bc}. 
The accuracy depends on  the domain size $x_\epsilon$ for $\Omega_\epsilon = [-x_\epsilon, x_\epsilon]$.

\begin{figure}[h!]
\hspace{-0.5cm}
\includegraphics[height=6cm]{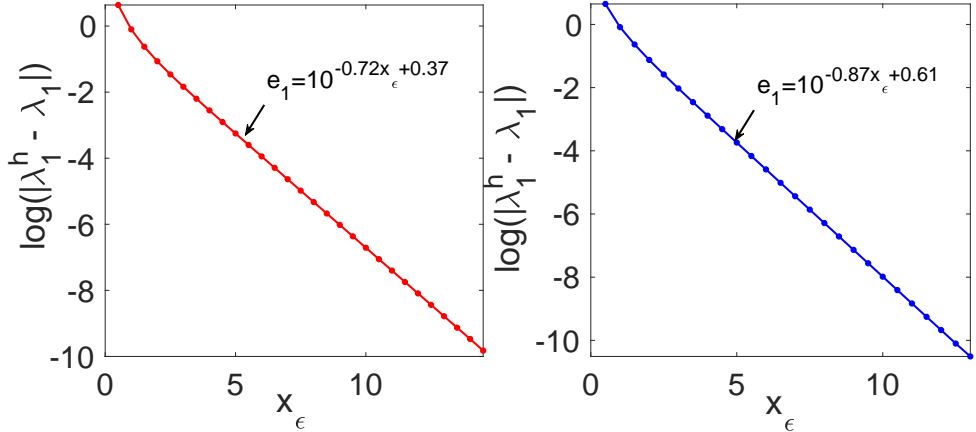}
\vspace{-0.8cm}
\caption{Eigenvalue errors versus domain size $x_\epsilon$ for the two-body problem with polynomially decaying potential \eqref{eq:pd} (left plot) and exponentially decaying potential \eqref{eq:ed} (right plot). }
\label{fig:2bpdom}
\end{figure}

\begin{figure}[h!]
\hspace{-0.5cm}
\includegraphics[height=6cm]{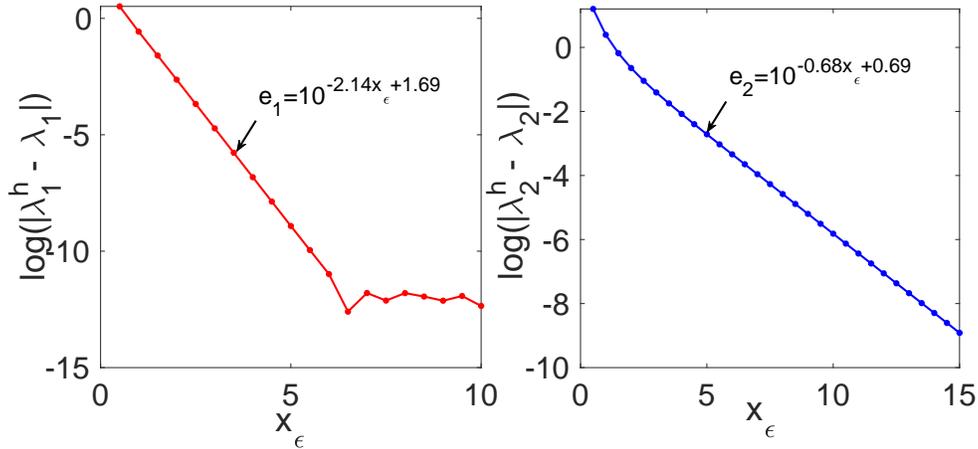}
\caption{Eigenvalue errors versus domain size $x_\epsilon$ for the two bound states of the two-body problem with exponentially decaying potential \eqref{eq:ed}.}
\label{fig:2bpdom2}
\end{figure}

To study the impact of the domain size $x_\epsilon$ on the accuracy,
we apply the high-accuracy IGA method with $C^6$ septic B-spline elements. 
We apply uniform mesh grids with a fixed grid size $h=0.01$.
This setting of using a high-order element with fine grid size is to guarantee that the errors are dominated by the choice of the domain size.
Figure \ref{fig:2bpdom} shows how the eigenvalue errors decrease when the domain size $x_\epsilon$ increases.
We set $\kappa=1/2$ in \eqref{eq:pde} and the potential magnitude $\beta=1$.
The left plot of Figure \ref{fig:2bpdom} shows that the eigenvalue error decays exponentially when $x_\epsilon$ increases for the potential \eqref{eq:pd} while the right plot of Figure \ref{fig:2bpdom} shows that of the potential \eqref{eq:ed}. 
The fitted functions that establish the relation between the error and the domain size are $e_1=10^{-0.72x_\epsilon + 0.37}$ and $e_1=10^{-0.87x_\epsilon + 0.61}$ for \eqref{eq:pd} and \eqref{eq:ed}, respectively. 
Figure \ref{fig:2bpdom2} shows the case where there are two bound states. Therein, potential \eqref{eq:ed} is used with $\beta=5$. We observe a similar behaviour. 
The errors of the ground state reach an order of $10^{-12}$ when $x_\epsilon \ge 6.5$. 
This is due to that the IGA discretization error dominates the overall error (from discretization and approximation of the domain). 
The fitted functions give guidance for choosing the domain $\Omega_\epsilon$ appropriately. 
For example, for the two-body problem with a potential \eqref{eq:ed} and $\beta=5$, the fitted function for the second bound state (a larger domain is required to compute this mode) is $e_2=10^{-0.68x_\epsilon + 0.69}$.
Thus, to achieve an accuracy of error $10^{-15}$, we set $10^{-15} = 10^{-0.68x_\epsilon + 0.69}$ and solve for $x_\epsilon$ to get the domain $\Omega_\epsilon = [-x_\epsilon, x_\epsilon] = [-23.1, 23.1]$.
This means that we require to solve the problem with a minimal domain size $x_\epsilon = 23.1$ to get an accuracy of order $10^{-15}$.

\vspace{-0.5cm}

\vspace{-0.5cm}

\subsection{Three-body problem}

Now, we consider the three-body problem with a heavy-light body ratio $m_h/m_l=20$ that is studied in \cite{happ2019universality,happ2022universality,thies2021exploiting}.
With such a mass ratio, $\kappa$ in the unified problem \eqref{eq:pde} is a matrix with entries $\kappa = (41/84, 0; 0, 2/21)$.
Our goal is to approximate the eigenenergies and eigenstates obtained by using the classical BO approximation in these papers. 
This preliminary numerical study demonstrates that the proposed IGA method is a promising alternative to the classical BO approximation method that a strong assumption is posed on the mass ratio. 

\begin{figure}[h!]
\includegraphics[height=10cm]{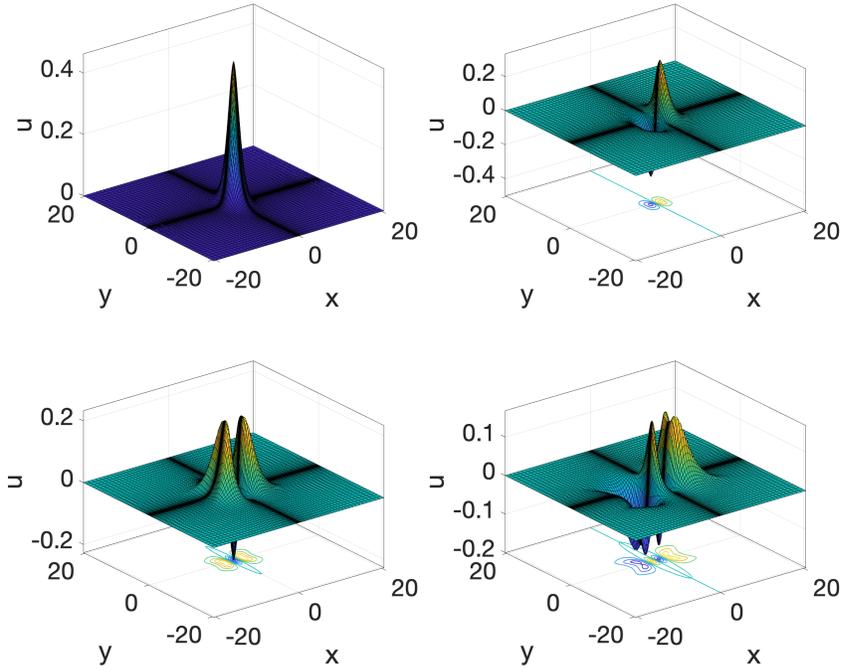}
\vspace{-0.8cm}
\caption{ Eigenstate wavefunctions of the three-body problem with a mass ratio 20, potential magnitude $\beta=0.344595351$, and the exponentially decaying potential \eqref{eq:ed}.}
\label{fig:3bp}
\end{figure}

To solve the three-body problem, we apply the highly accurate IGA method with $C^6$ septic B-spline elements. 
We set the domain as $\Omega_\epsilon = [-20, 20]$ and apply a non-uniform grid with $80\times80$ elements. 
Table \ref{tab:eig} shows the eigenvalues of the bound states of the three-body problem with the exponentially decaying potential \eqref{eq:ed}.
The potential magnitude is $\beta=0.344595351$. 
The IGA eigenvalues are close to the ones (scaled) shown in Table 1 of \cite{thies2021exploiting}. 
Figure \ref{fig:3bp} shows the first four bound state eigenfunctions.
The eigenstate solution shapes match well with the ones obtained using the BO approximation in Figure 4 of \cite{happ2019universality}. 
Moreover, we observe a similar universality behaviour as in \cite{happ2019universality} and we will present a detailed study in future work. 
In conclusion, the IGA method with a small mesh grid has the ability to approximate well both eigenenergies and eigenfunctions of the bound states of the three-body problem. 

\vspace{-0.3cm}
\begin{table}[ht]
\centering 
\begin{tabular}{|c|c| cc | cc |}
\hline
method & $\beta$ & $j$ (Bosons)  & $\lambda_j^h$ &  $j$ (Fermions) & $\lambda_j^h$   \\[0.1cm] \hline
&                       & 0 & -0.2476034576 & 1 & -0.1825896533 \\[0.1cm]
IGA &0.344595351  & 2 & -0.1412793292 & 3 & -0.1182591543 \\[0.1cm]
&                       & 4 & -0.1060931444 & 5 & -0.1005294105 \\[0.1cm] \hline

&                       & 0 & -0.247603458 & 1 & -0.182589653 \\[0.1cm]
BO in \cite{thies2021exploiting} &0.34459535  & 2 & -0.141279329 & 3 & -0.118259157 \\[0.1cm]
&                       & 4 & -0.106093864 & 5 & -0.102845702 \\[0.1cm] \hline
                       
%
\end{tabular}
\vspace{0.3cm}
\caption{Eigenenergies of the three-body problem with a mass ratio 20 and an exponentially decaying potential \eqref{eq:ed} when using IGA and BO approximation in \cite{thies2021exploiting}.}
\label{tab:eig} 
\end{table}

\vspace{-1cm}

\section{Concluding remarks} \label{sec:conclusion} 

In this paper, we initiated IGA of the quantum two- and three-body problems.  IGA is developed based on the classical Galerkin FEA 
that has the advantages of a mature theoretical understanding of the error estimates, stabilities, and robustness.
IGA is successfully applied to solve the bound states of the two- and three-body problems in 1D with arbitrary mass ratios and potential shapes.

As for future work, the first possible direction is a generalization to the two- and three-body problems in multiple dimensions. 
Tensor-product structures may be applied to reduce the computational costs \cite{thies2021exploiting}. 
Another direction of future work is that one may use the recently developed softFEM \cite{deng2021softfem} and dispersion-minimized blending quadratures \cite{deng2018dispersion,calo2019dispersion} to solve the three-body problem with higher accuracy and efficiency.


%
%
%
\bibliographystyle{splncs04}
\bibliography{ref.bib}
%


%
%
%
%

\end{document}